\documentclass[12pt]{article}

\usepackage{amssymb,amsbsy,amsmath,amsfonts,amssymb,amscd}
\usepackage{latexsym,euscript,exscale,epic,eepic,epsfig}
\usepackage{latexsym,euscript,exscale,epic,eepic,epsfig}
\usepackage[english,francais]{babel}
\newtheorem{theorem}{Theorem}[section]
\newtheorem{lemma}[theorem]{Lemma}
\newtheorem{proposition}[theorem]{Proposition}
\newtheorem{corollary}[theorem]{Corollary}
\newtheorem{e-definition}[theorem]{Definition\rm}

\def\og{\leavevmode\raise.3ex\hbox{$\scriptscriptstyle\langle\!\langle$~}}
\def\fg{\leavevmode\raise.3ex\hbox{~$\!\scriptscriptstyle\,\rangle\!\rangle$}}
\begin{document}
%\begin{frontmatter}
\selectlanguage{english}
\vskip 0.5\baselineskip
\title{Boundedness commutator between Sobolev spaces}
\author{By Sadek Gala\\Universit\'e d'Evry Val d'Essonne\\
  D\'epartement de math\'ematiques\\ Bd F. Mitterrand. 91025 Evry
  Cedex. France\\ Sadek.Gala@maths.univ-evry.fr }
\maketitle
%\ead{Sadek.Gala@maths.univ-evry.fr}
%\address{Universit\'e d'Evry Val d'Essonne\\ D\'epartement de math\'ematiques\\ Bd F. Mitterrand. 91025 Evry Cedex. France}

\begin{abstract} In this paper,we will study the boundedness properties of commutator
\[
C_{f}=\left[  f,\Delta\right]
\]
acting from $\overset{.}{H}^{1}\left(  \mathbb{R}^{d}\right)  $ to
$\overset{.}{H}^{-1}\left(  \mathbb{R}^{d}\right)  .$
%\end{frontmatter}
%\begin{frontmatter}
%\selectlanguage{francais}
\vskip 0.5\baselineskip
\noindent% Text of abstract in French

%\noindent{\bf R\'esum\'e}\\
%\begin{abstract}

%\selectlanguage{francais}
% Texte du rÈsumÈ en franÁais
%{\bf R\'esum\'e.}
%Dans cet article, on \'{e}tudiera la propri\'{e}t\'{e} de
%continuit\'{e} du commutateur agissant de l'espace de Sobolev $\overset{.}%
%{H}^{1}\left(  \mathbb{R}^{d}\right)  $ sur son dual $\overset{.}{H}%
%^{-1}\left(  \mathbb{R}^{d}\right)  .$
\end{abstract}
%\end{frontmatter}

%\section{Introduction}
Throughout this paper, $d$ will denote a fixed positive integer greater than 2
and $\Omega$ will denote an open, non-empty subset of $\mathbb{R}^{d}$. We
wish to study linear elliptic differential operator $\Delta$ in divergence
form defined on $C^{2}$ functions by
\[
\Delta f(x)=\frac{1}{2}\sum_{i,j=1}^{d}\partial_{i}\left(  \partial
_{j}f\right)  (x)
\]

Recall the divergence theorem. Suppose $\Omega$ is a nice region, $F$ is a
smooth vector field, $\nu(x)$ is the outward normal vector at $x\in
\partial\Omega$, and $\sigma$ is surface measure on $\partial\Omega$. The
divergence theorem then says that%
\[
\int_{\partial\Omega}F(y).\nu(y)\sigma(dy)=\int_{\Omega}divF(x)dx.
\]

\bigskip A twice continuously differentiable, complex-valued function $v$
defined on $\Omega$ is harmonic on $\Omega$ if
\[
\Delta v=0.
\]

We begin with the following simple observation :

\begin{proposition}
\label{prop1}
Let $g$ be a $C^{\infty}$ function with compact support and $f$ a bounded
$C^{\infty}$ function. Then
\[
\int_{\mathbb{R}^{d}}g(x)\Delta f(x)dx=-\frac{1}{2}\int_{\mathbb{R}^{d}%
}\left(  \sum_{i,j=1}^{d}\partial_{i}g(x)\partial_{j}f(x)\right)  dx.
\]
The integrand on the right could be written $\nabla g.\nabla f.$
\end{proposition}

\begin{proof}
We apply the divergence theorem. Let $B$ be a ball large enough to contain the
support of $g$ and let $F(x)$ be the vector field whose $i-$th component is
\[
\frac{1}{2}g(x)\sum_{j=1}^{d}\partial_{j}f(x).
\]
Since $g$ is $0$ on $\partial B$, then $F.\nu=0$ on $\partial B$, and also,
\begin{align*}
divF(x) &  =\frac{1}{2}\sum_{i=1}^{d}\partial_{i}\left(  g(x)\sum_{j=1}%
^{d}\partial_{j}f(x)\right) =\frac{1}{2}\sum_{i,j=1}^{d}\partial_{i}g(x)\partial_{j}f(x)+g(x)\Delta
f(x).
\end{align*}
\end{proof}

\begin{proposition}
\label{prop2}
Suppose $v$ is positive and harmonic in $Q\left(  4\right)  $.
There exists $c_{1}$ independent of $v$ such that if $w=\log v$, then
\[
\int_{Q}\left|  \nabla w(x)\right|  ^{2}dx\leq c_{1}h^{d-2}%
\]
for all cubes $Q$ of side length $h$ contained in $Q(2)$.
\end{proposition}

\begin{proof}
Let $Q^{\ast}$ be the cube with the same center as $Q$ but side length twice
as long. Note $Q^{\ast}\subset Q\left(  4\right)  $. Let $\varphi$ be
$C^{\infty}$ with values in $\left[  0,1\right]  $, equal to $1$ on $Q$,
supported in $Q^{\ast}$, and such that $\left\|  \nabla\varphi\right\|
_{\infty}\leq c_{2}h^{-1}$. Since
\[
\nabla w=\frac{\nabla v}{v}%
\]
and $v$ is harmonic in $Q\left(  4\right)  $,%
\begin{align*}
0  =2\int\frac{\varphi^{2}}{v}\Delta v=-\int\nabla\left(  \frac{\varphi
^{2}}{v}\right)  .\nabla v =-\int2\frac{\varphi\nabla\varphi}{v}.\nabla v+\int\frac{\varphi^{2}}%
{v^{2}}\nabla v.\nabla v \\=-2\int\varphi\nabla\varphi.\nabla w+\int\varphi^{2}\nabla w.\nabla w.
\end{align*}
So by the Cauchy-Schwarz inequality and proposition \ref{prop1},%
\begin{align*}
\int_{Q^{\ast}}\varphi^{2}\left|  \nabla w\right|  ^{2}dx  &  \leq c_{3}%
\int_{Q^{\ast}}\varphi^{2}\nabla w.\nabla w=c_{4}\int_{Q^{\ast}}\nabla
\varphi.\varphi\nabla w\\
&  \leq c_{5}\left(  \int_{Q^{\ast}}\left|  \nabla\varphi\right|
^{2}dx\right)  ^{\frac{1}{2}}\left(  \int_{Q^{\ast}}\varphi^{2}\left|  \nabla
w\right|  ^{2}dx\right)  ^{\frac{1}{2}}.
\end{align*}
Dividing by the second factor on the right, squaring and using the bound on
$\left|  \nabla\varphi\right|  $,%
\[
\int_{Q}\left|  \nabla w(x)\right|  ^{2}dx\leq\left(  \int_{Q^{\ast}}%
\varphi^{2}\left|  \nabla w\right|  ^{2}dx\right)  \leq c_{5}^{2}\left|
Q^{\ast}\right|  c_{2}^{2}h^{-2},
\]
which implies our result.
\end{proof}

We need to distinguish the class of vector fields $\overrightarrow{F}$ such
that the commutator inequality
\[
\left|  \int_{\mathbb{R}^{d}}\overrightarrow{F}.\left(  \overline{u}\nabla
v-v\nabla\overline{u}\right)  dx\right|  \leq C\left\|  u\right\|
_{\overset{.}{H}^{1}}\left\|  v\right\|  _{\overset{.}{H}^{1}}%
\]
for all $u$, $v\in\mathcal{D}\left(  \mathbb{R}^{d}\right)  $. In the
important case where $\overrightarrow{F}=\nabla f$, the preceding inequality
is equivalent to the boundedness of the commutator
\[
C_{f}=\left[  f,\Delta\right]
\]
acting from $\overset{.}{H}^{1}\left(  \mathbb{R}^{d}\right)  $ to
$\overset{.}{H}^{-1}\left(  \mathbb{R}^{d}\right)  .$

Another useful point in the sequel is the symmetry of the Laplace operator,
\[
\left\langle \Delta\varphi,\psi\right\rangle =-\left\langle D\varphi
,D\psi\right\rangle =\left\langle \varphi,\Delta\psi\right\rangle
\]
for all $\varphi$, $\psi\in\mathcal{D}\left(  \mathbb{R}^{d}\right)  $.

Let us present a proof of the 

\begin{lemma}
For any twice continuously differentiable functions $u:A\rightarrow\mathbb{R}$
and $v:A\rightarrow\mathbb{R}$, where $A\subset\mathbb{R}^{d}$, we have%

\[
div\left(  u\nabla v-v\nabla u\right)  =u\Delta v-v\Delta u.
\]
\end{lemma}

\begin{proof}
Note that
\begin{align*}
&  div\left(  u\nabla v-v\nabla u\right) =div\left(  u\left(  \frac{\partial v}{\partial x_{1}},...,\frac{\partial
v}{\partial x_{d}}\right)  -v\left(  \frac{\partial u}{\partial x_{1}%
},...,\frac{\partial u}{\partial x_{d}}\right)  \right) \\
&  =div\left(  u\frac{\partial v}{\partial x_{1}}-v\frac{\partial u}{\partial
x_{1}},...,u\frac{\partial v}{\partial x_{d}}-v\frac{\partial u}{\partial
x_{d}}\right) \\
&  =\frac{\partial}{\partial x_{1}}\left(  u\frac{\partial v}{\partial x_{1}%
}-v\frac{\partial u}{\partial x_{1}}\right)  +...+\frac{\partial}{\partial
x_{d}}\left(  u\frac{\partial v}{\partial x_{d}}-v\frac{\partial u}{\partial
x_{d}}\right) \\
&  =\left(  u\frac{\partial^{2}v}{\partial x_{1}^{2}}+\frac{\partial
u}{\partial x_{1}}\frac{\partial v}{\partial x_{1}}-\frac{\partial v}{\partial
x_{1}}\frac{\partial u}{\partial x_{1}}-v\frac{\partial^{2}u}{\partial
x_{1}^{2}}\right)  +...+\left(  u\frac{\partial^{2}v}{\partial x_{d}^{2}%
}+\frac{\partial u}{\partial x_{d}}\frac{\partial v}{\partial x_{d}}%
-\frac{\partial v}{\partial x_{d}}\frac{\partial u}{\partial x_{d}}%
-v\frac{\partial^{2}u}{\partial x_{d}^{2}}\right) \\
&  =u\left(  \frac{\partial^{2}v}{\partial x_{1}^{2}}+...+\frac{\partial^{2}%
v}{\partial x_{d}^{2}}\right)  -v\left(  \frac{\partial^{2}u}{\partial
x_{1}^{2}}+...+\frac{\partial^{2}u}{\partial x_{d}^{2}}\right) =u\Delta v-v\Delta u.
\end{align*}
\end{proof}

\begin{lemma}
Let $f\in\mathcal{D}^{\prime}\left(  \mathbb{R}^{d}\right)  $. Then for all
$u$, $v\in\mathcal{D}\left(  \mathbb{R}^{d}\right)  $, we have
\[
\left|  \left\langle C_{f}u,v\right\rangle\right|  =\left|  \int
_{\mathbb{R}^{d}}\overrightarrow{F}.\left(  \overline{u}\nabla v-v\nabla
\overline{u}\right)  dx\right|  .
\]
\end{lemma}

\begin{proof}
Let $u$, $v\in\mathcal{D}\left(  \mathbb{R}^{d}\right)  .$ We observe that%
\begin{align*}
\left\langle C_{f}u,v\right\rangle  &  =\left\langle \left[  f,\Delta\right]
u,v\right\rangle =\left\langle f\Delta u-\Delta\left(  fu\right)
,v\right\rangle \\
&  =\left\langle f\Delta u,v\right\rangle -\left\langle \Delta\left(
fu\right)  ,v\right\rangle =\left\langle f,v\Delta\overline{u}\right\rangle
-\left\langle fu,\Delta v\right\rangle \\
&  =\left\langle f,v\Delta\overline{u}\right\rangle -\left\langle
f,\overline{u}\Delta v\right\rangle =\left\langle f,v\Delta\overline
{u}-\overline{u}\Delta v\right\rangle \\
&  =\left\langle f,div\left(  v\nabla\overline{u}-\overline{u}\nabla v\right)
\right\rangle \\
&  =-\left\langle \nabla f,v\nabla\overline{u}-\overline{u}\nabla
v\right\rangle =-\left\langle \overrightarrow{F},v\nabla\overline{u}%
-\overline{u}\nabla v\right\rangle
\end{align*}
Therefore%
\[
\left|  \left\langle C_{f}u,v\right\rangle \right|  =\left|  \int
_{\mathbb{R}^{d}}\overrightarrow{F}.\left(  \overline{u}\nabla v-v\nabla
\overline{u}\right)  dx\right|  .
\]
\end{proof}

We will let $P$ be the Newton potential of $\mu$ defined by%
\[
P\mu(x)=P(x)=C(d)\int_{\mathbb{R}^{d}}\frac{d\mu(y)}{\left|  x-y\right|
^{d-2}},\text{ \ \ }x\in\mathbb{R}^{d}\text{, \ \ }d\geq3.
\]
Then, we have the following lemma :

\begin{lemma}
\label{lemma 2.7}Let $\mu$ be a positive Borel measure on $\mathbb{R}^{d}$
such that $P(x)=I_{2}\mu(x)\neq\infty$. Then the following inequality hold :%

\begin{equation}
\int_{\mathbb{R}^{d}}u^{2}(x)\frac{\left|  \nabla P(x)\right|  ^{2}}{P^{2}%
(x)}dx\leq4\left\|  \nabla u\right\|  _{L^{2}}^{2}\text{, \ \ }u\in
\mathcal{D}\left(  \mathbb{R}^{d}\right)  .\label{eq2}%
\end{equation}
\end{lemma}

\begin{proof}
Suppose $u\in\mathcal{D}\left(  \mathbb{R}^{d}\right)  .$ Then $K=$supp $u$ is
a compact set, and obviously $\underset{x\in K}{\inf}P(x)>0.$ Without loss of
generality, we assume that $\nabla P\in L_{loc}^{2}\left(  \mathbb{R}%
^{d}\right)  $, and hence the left-hand side of (\ref{eq2}) is finite.

Using integration by parts, together with the properties $-\Delta P=\mu$
(understood in the distributional sense) and applying the Schwarz inequality,
we get :%
\begin{align*}
2\int_{\mathbb{R}^{d}}\nabla u(x).\nabla P(x)\frac{u(x)}{P(x)}dx  &
=\int_{\mathbb{R}^{d}}u^{2}(x)\frac{\left|  \nabla P(x)\right|  ^{2}}%
{P^{2}(x)}dx+\int_{\mathbb{R}^{d}}u^{2}(x)\frac{d\mu(x)}{P(x)}dx\\
&  \leq2\left(  \int_{\mathbb{R}^{d}}u^{2}(x)\frac{\left|  \nabla P(x)\right|
^{2}}{P^{2}(x)}dx\right)  ^{\frac{1}{2}}\left(  \int_{\mathbb{R}^{d}}\left|
\nabla u(x)\right|  ^{2}dx\right)  ^{\frac{1}{2}},
\end{align*}
for all $u\in\mathcal{D}\left(  \mathbb{R}^{d}\right)  .$ Using this
inequality, we have
\[
\int_{\mathbb{R}^{d}}u^{2}(x)\frac{\left|  \nabla P(x)\right|  ^{2}}{P^{2}%
(x)}dx\leq2\left(  \int_{\mathbb{R}^{d}}u^{2}(x)\frac{\left|  \nabla
P(x)\right|  ^{2}}{P^{2}(x)}dx\right)  ^{\frac{1}{2}}\left(  \int
_{\mathbb{R}^{d}}\left|  \nabla u(x)\right|  ^{2}dx\right)  ^{\frac{1}{2}}.
\]
Consequently, we obtain%
\[
\int_{\mathbb{R}^{d}}u^{2}(x)\frac{\left|  \nabla P(x)\right|  ^{2}}{P^{2}%
(x)}dx\leq4\int_{\mathbb{R}^{d}}\left|  \nabla u(x)\right|  ^{2}dx.
\]
\end{proof}

\begin{lemma}
\bigskip\label{lemma 2}The Sobolev space $\overset{.}{H}^{1}\left(
\mathbb{R}^{d}\right)  $ be invariant under the transformation :%

\begin{align*}
u  & \rightarrow\widetilde{u}=e^{iw}u\text{, \ \ \ }v\rightarrow\widetilde
{v}=e^{iw}v,\\
\left\|  u\right\|  _{\overset{.}{H}^{1}}  & \simeq\left\|  \widetilde
{u}\right\|  _{\overset{.}{H}^{1}},\text{ \ \ \ \ \ \ }\left\|  v\right\|
_{\overset{.}{H}^{1}}\simeq\left\|  \widetilde{v}\right\|  _{\overset{.}%
{H}^{1}}%
\end{align*}
\end{lemma}

\begin{proof}
Clearly, We observe that $v\in\overset{.}{H}^{1}\left(  \mathbb{R}^{d}\right)
$
\[
\nabla\left(  e^{iw}u\right)  =\left(  iu\nabla w+\nabla u\right)  e^{iw}.
\]
Consequently, for every $u\in\overset{.}{H}^{1}\left(  \mathbb{R}^{d}\right)
,$%
\[
\left\|  \widetilde{u}\right\|  _{\overset{.}{H}^{1}}=\left\|  e^{iw}%
u\right\|  _{\overset{.}{H}^{1}}\leq\left\|  u\nabla w\right\|  _{L^{2}%
}+\left\|  \nabla u\right\|  _{L^{2}}.
\]
We set 
\[
w=\log P\mu\text{, \ \ }v=P\mu.
\]
Note that
\[
\nabla w=\frac{\nabla P}{P}.
\]
Hence, by lemma \ref{lemma 2.7},
\[
\left\|  u\nabla w\right\|  _{\overset{.}{H}^{1}}=\left\|  u\frac{\nabla P}%
{P}\right\|  _{\overset{.}{H}^{1}}\leq2\left\|  u\right\|  _{\overset{.}%
{H}^{1}}.
\]
From this, it follows
\[
\left\|  \widetilde{u}\right\|  _{\overset{.}{H}^{1}}\leq3\left\|  u\right\|
_{\overset{.}{H}^{1}}.
\]
Using similar estimates for $\widetilde{u}\rightarrow u=e^{-iw}\widetilde{u}$,
we deduce
\[
\frac{1}{3}\left\|  u\right\|  _{\overset{.}{H}^{1}}\leq\left\|  \widetilde
{u}\right\|  _{\overset{.}{H}^{1}}.
\]
This completes the proof of lemma.
\end{proof}

We now establish the main result of this paper.

\begin{theorem}
\label{theorem}
Let $\overrightarrow{F}\in\mathcal{D}^{\prime}\left(  \mathbb{R}^{d}\right)
$, $d\geq3.$ Suppose that, for all $u$, $v\in\mathcal{D}\left(  \mathbb{R}%
^{d}\right)  $%

\begin{equation}
\left|  \left\langle\overrightarrow{F},v\nabla\overline{u}-\overline{u}\nabla
v\right\rangle\right|  \leq C\left\|  u\right\|  _{\overset{.}{H}^{1}%
}\left\|  v\right\|  _{\overset{.}{H}^{1}}.\label{eq4}%
\end{equation}
Then there is a vector-field $\overrightarrow{G}=\nabla\left(  \Delta
^{-1}div\overrightarrow{F}\right)  \in L_{loc}^{2}\left(  \mathbb{R}%
^{d}\right)  $ such that
\[
\int_{\mathbb{R}^{d}}\left|  \overrightarrow{G}(x)\right|  ^{2}\left|
u(x)\right|  ^{2}dx\leq C\left\|  u\right\|  _{\overset{.}{H}^{1}}^{2}%
\]
for all $u\in\mathcal{D}\left(  \mathbb{R}^{d}\right)  $, where the constant
$C$ does not depend on  $u\in\mathcal{D}\left(  \mathbb{R}^{d}\right)  .$
\end{theorem}

\begin{proof}
Suppose that
\[
\left|  \left\langle \overrightarrow{F},v\nabla\overline{u}-\overline{u}\nabla
v\right\rangle \right|  \leq C\left\|  u\right\|  _{\overset{.}{H}^{1}%
}\left\|  v\right\|  _{\overset{.}{H}^{1}}.
\]
holds. By the continuity the bilinear form on the left-hand side can be
extended to all $u$, $v\in\overset{.}{H}^{1}\left(  \mathbb{R}^{d}\right)  $.
Let $v$ be a nonnegative function such that $w=\log$ $v$ and by the
proposition \ref{prop1}
\[
\nabla w=\frac{\nabla v}{v}\in L_{loc}^{2}\left(  \mathbb{R}^{d}\right)  .
\]
We need $w$ to be chosen so that the Sobolev space $\overset{.}{H}^{1}\left(
\mathbb{R}^{d}\right)  $ be invariant under the transformation :%
\begin{align*}
u &  \rightarrow\widetilde{u}=e^{iw}u\text{, \ \ \ }v\rightarrow\widetilde
{v}=e^{iw}v,\\
\left\|  u\right\|  _{\overset{.}{H}^{1}} &  \simeq\left\|  \widetilde
{u}\right\|  _{\overset{.}{H}^{1}},\text{ \ \ \ \ \ \ }\left\|  v\right\|
_{\overset{.}{H}^{1}}\simeq\left\|  \widetilde{v}\right\|  _{\overset{.}%
{H}^{1}}%
\end{align*}
Applying lemma \ref{lemma 2} and (\ref{eq4}), with $\widetilde{u}$ and
$\widetilde{v}$ in place of $u$ and $v$, we get%
\begin{align*}
\left|  \left\langle \overrightarrow{F},\widetilde{v}\nabla\overline
{\widetilde{u}}-\overline{\widetilde{u}}\nabla\widetilde{v}\right\rangle
\right|   &  \leq C\left\|  \widetilde{u}\right\|  _{\overset{.}{H}^{1}%
}\left\|  \widetilde{v}\right\|  _{\overset{.}{H}^{1}}\\
&  \leq9C\left\|  u\right\|  _{\overset{.}{H}^{1}}\left\|  v\right\|
_{\overset{.}{H}^{1}}.
\end{align*}
Notice that
\[
\widetilde{v}\nabla\overline{\widetilde{u}}-\overline{\widetilde{u}}%
\nabla\widetilde{v}=v\nabla\overline{u}-\overline{u}\nabla v-2i\overline
{u}v\nabla w.
\]
Combining the preceding estimates, we obtain%
\begin{align*}
2\left|  \left\langle \overrightarrow{F},\overline{u}v\nabla w\right\rangle
\right|   &  \leq\left|  \left\langle \overrightarrow{F},v\nabla\overline
{u}-\overline{u}\nabla v\right\rangle \right|  +\left|  \left\langle
\overrightarrow{F},\widetilde{v}\nabla\overline{\widetilde{u}}-\overline
{\widetilde{u}}\nabla\widetilde{v}\right\rangle \right|  \\
&  \leq10C\left\|  u\right\|  _{\overset{.}{H}^{1}}\left\|  v\right\|
_{\overset{.}{H}^{1}}.
\end{align*}
Observe that
\[
v\nabla w=\nabla v.
\]
Thus, we arrive at the inequality%
\[
\left|  \left\langle \overrightarrow{F},\overline{u}\nabla v\right\rangle
\right|  \leq C^{\prime}\left\|  u\right\|  _{\overset{.}{H}^{1}}\left\|
v\right\|  _{\overset{.}{H}^{1}}.
\]
From the preceding estimate and (\ref{eq4}), we deduce%
\begin{align*}
\left|  \left\langle \overrightarrow{F},\nabla\left(  \overline{u}v\right)
\right\rangle \right|   &  =\left|  \left\langle \overrightarrow{F}%
,v\nabla\overline{u}+\overline{u}\nabla v\right\rangle \right|  \leq\left|
\left\langle \overrightarrow{F},v\nabla\overline{u}\right\rangle \right|
+\left|  \left\langle \overrightarrow{F},\overline{u}\nabla v\right\rangle
\right|  \\
&  \leq C"\left\|  u\right\|  _{\overset{.}{H}^{1}}\left\|  v\right\|
_{\overset{.}{H}^{1}}.
\end{align*}
This yields
\[
\left|  \left\langle \left(  div\overrightarrow{F}\right)  u,v\right\rangle
\right|  \leq C"\left\|  u\right\|  _{\overset{.}{H}^{1}}\left\|  v\right\|
_{\overset{.}{H}^{1}},
\]
where $u\in\overset{.}{H}^{1}\left(  \mathbb{R}^{d}\right)  $ and $v=P\mu$. By
theorem 2.2 in [MV],%
\[
div\overrightarrow{F}\in\mathcal{M}\left(  \overset{.}{H}^{1}\rightarrow
\overset{.}{H}^{-1}\right)  .
\]
Hence, there is a vector-field $\overrightarrow{G}\in L_{loc}^{2}\left(
\mathbb{R}^{d}\right)  $ such that
\[
\overrightarrow{G}=\nabla\left(  \Delta^{-1}div\overrightarrow{F}\right)
=\nabla f\in\mathcal{M}\left(  \overset{.}{H}^{1}\rightarrow L^{2}\right)  .
\]
\end{proof}

\begin{corollary}
Under the assumtions of theorem \ref{theorem}, it follows that if
\[
\left|  \left\langle \overrightarrow{F},v\nabla\overline{u}-\overline{u}\nabla
v\right\rangle \right|  \leq C\left\|  u\right\|  _{\overset{.}{H}^{1}%
}\left\|  v\right\|  _{\overset{.}{H}^{1}},
\]
then
\[
\Delta f\in\mathcal{M}\left(  \overset{.}{H}^{1}\rightarrow\overset{.}{H}%
^{-1}\right)  .
\]
\end{corollary}

\begin{theorem}
Let $f\in\mathcal{D}^{\prime}\left(  \mathbb{R}^{d}\right)  .$ Then the
inequality
\[
\left|  \left\langle \left[  f,\Delta\right]  u,v\right\rangle \right|  \leq
C\left\|  u\right\|  _{\overset{.}{H}^{1}}\left\|  v\right\|  _{\overset{.}%
{H}^{1}},
\]
holds for all $u,v\in\mathcal{D}\left(  \mathbb{R}^{d}\right)  $, if and only
if there is a vector-field $\overrightarrow{G}\in L_{loc}^{2}\left(
\mathbb{R}^{d}\right)  $ such that%
\[
\int_{\mathbb{R}^{d}}\left|  \overrightarrow{G}(x)\right|  ^{2}\left|
u(x)\right|  ^{2}dx\leq C\left\|  u\right\|  _{\overset{.}{H}^{1}}^{2},
\]
where $\overrightarrow{G}=\nabla f.$
\end{theorem}

\begin{proof}
Suppose that $\overrightarrow{G}=\nabla\left(  \Delta^{-1}div\overrightarrow
{F}\right)  $ and
\[
\int_{\mathbb{R}^{d}}\left|  \overrightarrow{G}(x)\right|  ^{2}\left|
u(x)\right|  ^{2}dx\leq C\left\|  u\right\|  _{\overset{.}{H}^{1}}^{2}.
\]
Then $\overrightarrow{G}.\nabla$ is form bounded, since by Schwarz's
inequality, we have
\[
\left|  \left\langle \overrightarrow{G}.\nabla u,v\right\rangle \right|
\leq\left\|  u\overrightarrow{G}\right\|  _{L^{2}}\left\|  \nabla v\right\|
_{L^{2}}\leq C\left\|  u\right\|  _{\overset{.}{H}^{1}}\left\|  v\right\|
_{\overset{.}{H}^{1}}.
\]
Hence, this inequality obviously yields
\begin{align*}
\left|  \left\langle \overrightarrow{F},v\nabla\overline{u}-\overline{u}\nabla
v\right\rangle \right|    & \leq\left|  \left\langle \overrightarrow
{F},v\nabla\overline{u}\right\rangle \right|  +\left|  \left\langle
\overrightarrow{F},\overline{u}\nabla v\right\rangle \right|  \\
& \leq C\left\|  u\right\|  _{\overset{.}{H}^{1}}\left\|  v\right\|
_{\overset{.}{H}^{1}}.
\end{align*}
This completes the proof of theorem.
\end{proof}

\end{document}